\documentclass[12pt]{article}
\usepackage{amsmath,amsthm,amssymb,amscd}
\usepackage{latexsym}
\newtheorem{thm}{Theorem}[section]
 \newtheorem{cor}[thm]{Corollary}
 
 \newtheorem{prop}[thm]{Proposition}
 \theoremstyle{definition}
 
 \theoremstyle{remark}

 \numberwithin{equation}{section}

\begin{document}

\title
{Local derivative estimates for the heat equation coupled to the Ricci flow}

\author{ Hong Huang}
\date{}
\maketitle
\begin{abstract}
  In this note we obtain local derivative estimates of Shi-type for the heat equation coupled to the Ricci flow. As applications,  in part combining with Kuang's work, we extend some results of Zhang and Bamler-Zhang including distance distortion estimates and a backward pseudolocality theorem for Ricci flow on compact manifolds to the noncompact case.



\end{abstract}

% ----------------------------------------------------------------

\section {Introduction}

The Bernstein method is a strategy for obtaining derivative estimates for a solution to some PDE via applying the maximum principle to some partial differential inequality satisfying by a suitable combination of the solution and its derivatives. It is very useful in PDE and geometric analysis. In \cite{B87}/\cite{S89} Bando/Shi adapted this method to obtain global/local derivative estimates for the Ricci flow. Shi's local derivative estimates are fundamental for the study of the Ricci flow. For expositions and/or alternative proofs of Shi's local derivative estimates see for example, Hamilton \cite{H95}, Cao-Zhu \cite{CZ}, Chow-Lu-Ni \cite{CLN}, Chow et al \cite{C+08} and Tao \cite{T}.
With bounds on some derivatives of curvatures of the initial metrics Lu (see \cite{LT}, \cite{MT} and \cite{C+08}) got a  modified version of Shi's local derivative estimates.  Ecker-Huisken \cite{EH} got  Shi-type estimates for the mean curvature flow.
Grayson-Hamilton \cite{GH} derived Shi-type estimates  for  the harmonic map heat flow. For the heat equation on a Riemannian manifold, Kotschwar \cite{K} obtained  a Shi-type local gradient estimate,  while  the author \cite{Hu} obtained  local higher order derivative estimates. Recently the Shi-type estimates are  also derived for some other geometric evolution equations. See for example \cite{LW} and \cite{Ch}.

In his lectures at Tsinghua University in 2012/13 and in 2015, Hamilton \cite{H13} constructed a  comparison function (see Corollary 5.2 here), and used it to simplify  Shi's proof of the local derivative estimates for the Ricci flow. (Actually Hamilton \cite{H13} only gave the details for the gradient estimate. For  completeness here we  write down some  details for  the higher  derivative case by using Hamilton's comparison function; see the proof of Theorem 5.3.)

In this note we use Hamilton's comparison function to obtain local derivative estimates of Shi-type for the heat equation coupled with the Ricci flow. To state our results we first  introduce some notations.
 Fix $T >0$. Let $(M, (g(t))_{t\in [0,T]})$ be a solution (not necessarily complete) to Hamilton's Ricci flow
 \begin{equation*}
 \frac{\partial g(t)}{\partial t}=-2Ric (g(t))
 \end{equation*}
 on a manifold $M$ (without boundary) of dimension $n$.  For $x\in M$, $t\in [0,T]$ and $r>0$, let $B(x, t, r)$ be the open metric ball with center $x$ and of radius $r$ w.r.t. the metric $g(t)$,  and let $PB_r(x,T)$ denote a parabolic cylinder (as called in \cite{H13}) that is the set of all points $(x', t)$ with $ x' \in \overline{B(x,t,r)}$ (the closure of $B(x,t,r)$) and $t \in [0, T]$.  For any points $x, y \in M$ let $d_t(x, y)$ be the distance between   $x$ and $y$  w.r.t. $g(t)$.

 We have the following gradient estimate.
\begin{thm} \label{thm 1.1} \   Fix $T >0$. Let $(M, (g(t))_{t\in [0,T]})$ be a solution (not necessarily complete) to the Ricci flow on a manifold $M$ (without boundary) of dimension $n$.  Fix $x_0 \in M$ and $r>0$. Assume that the parabolic cylinder $PB_r(x_0,T)$ is compact, and $Ric\leq \frac{n-1}{r^2}$ on  $PB_r(x_0, T)$. Let $u$ be a  smooth solution to the heat equation $(\frac{\partial}{\partial t}-\Delta_{g(t)})u=0$ coupled to the Ricci flow on $M \times [0,T]$. Suppose $|u|\leq a$ on $PB_r(x_0, T)$, where $a$ is a positive constant. Then
\begin{equation*}
|\nabla u| \leq C_1a(\frac{1}{r}+\frac{1}{\sqrt{t}})    \hspace{2mm}  \text{on}   \hspace{2mm}    PB_{ \frac{r}{2}}(x_0, T)\setminus \{(x, 0) | x \in M\},
\end{equation*}
where the constant $C_1$ depends only on the dimension.
\end{thm}

Note that in  Theorem 2.2 in Bailesteanu-Cao-Pulemotov \cite{BCP}, where a two-sided bound on the Ricci curvature in a parabolic cylinder is assumed, a gradient estimate for a positive  solution $u$ to the heat equation coupled with the Ricci flow is given in terms of the pointwise value of $u$ and an upper bound of $u$ in the parabolic cylinder.  With a two-sided bound on the Ricci curvature of the form $|Ric|\leq \frac{n-1}{r^2}$ on  $PB_r(x_0, T)$ and given a positive  solution $u$ to the heat equation coupled with the Ricci flow, at a point where the value of  $u$ is very small,  the estimate in  Theorem 2.2 in \cite{BCP} is sharper than ours, but at a point where the value of  $u$ is not so small, the two estimates are comparable. See also Exercise 2.19 in Chow-Lu-Ni \cite{CLN} for a related global estimate.

We also get a Hessian estimate.
\begin{thm} \label{thm 1.2}   Let $M$ be  a manifold  (without boundary) of dimension $n$. Suppose $g(t)$ is a solution (not necessarily complete) to the Ricci flow  on $M \times [0, T]$ for some $T >0$.  Fix $x_0 \in M$ and $r>0$.  Assume that  the parabolic cylinder $PB_r(x_0,T)$ is compact, and $|Rm|\leq \frac{1}{r^2}$ on  $PB_r(x_0,  T)$. Let $ u $ be a smooth solution to the heat equation $(\frac{\partial}{\partial t}-\Delta_{g(t)})u=0$ coupled to the Ricci flow on $M \times [0, T]$. Suppose $|u|\leq a$  on $PB_r(x_0, T)$, where $a$ is a positive constant. Then
\begin{equation*}
|\nabla^2 u| \leq C_2a(\frac{1}{r^{2}}+\frac{1}{t})    \hspace{2mm}  \text{on}   \hspace{2mm}    PB_\frac{r}{4}(x_0,  T)\setminus \{(x, 0) | x \in M\},
\end{equation*}
where the constant $C_2$ depends only on the dimension.
\end{thm}
Compare Theorem 1.3 (b) in Han-Zhang \cite{HZ}, where an upper bound for the Hessian matrix of $u$ is obtained at points with certain distances away from the parabolic boundary.  One can also find a  global Hessian estimate in the proof of Theorem 18.2 in \cite{H95} assuming an initial gradient bound.

In Section 2 we prove Theorems 1.1 and 1.2, and derive similar estimates  for higher  derivatives.  In Section 3 we extend some derivative estimates in Zhang \cite{Z06}, Cao-Hamilton \cite{CH} and Bamler-Zhang \cite{BZ} on compact manifolds to the noncompact case using Theorems 1.1 and 1.2.  In Section 4 we get a slight improvement of some results in Kuang \cite{Ku1} and \cite{Ku2}, and point out that combining this and results in Section 3 one can extend some results in Zhang \cite{Z12a} and Bamler-Zhang \cite{BZ} including distance distortion estimates and a backward pseudolocality theorem for Ricci flow on compact manifolds to certain  noncompact manifolds.  Finally, in Section 5, which is an appendix,  we recall Hamilton's construction of a comparison function and the application to (the first order derivative case of) Shi's local derivative estimates, which appear in   \cite{H13}, and add some details for the higher order derivative case.

\vspace*{0.4cm}

\section{ Shi-type estimates}

 Fix $T >0$. Let $(M, (g(t))_{t\in [0,T]})$ be a solution (not necessarily complete) to the Ricci flow
 on a manifold $M$ (without boundary) of dimension $n$.
Let $u$ be a  smooth solution to the heat equation $(\frac{\partial}{\partial t}-\Delta_{g(t)})u=0$ coupled to the Ricci flow.  Recall that (compare for example \cite{CLN}, \cite{H95} and \cite{To})

\begin{equation}
(\frac{\partial}{\partial t}-\Delta)|\nabla u|^2=-2|\nabla^2u|^2,
\end{equation}

\begin{equation*}
(\frac{\partial}{\partial t}-\Delta)\nabla^2u=Rm*\nabla^2u,
\end{equation*}

\begin{equation}
(\frac{\partial}{\partial t}-\Delta)|\nabla^2u|^2=-2|\nabla^3u|^2+Rm*\nabla^2u*\nabla^2u,
\end{equation}

\begin{equation*}
(\frac{\partial}{\partial t}-\Delta)\nabla^ku=\sum_{i=0}^{k-2} \nabla^i Rm * \nabla^{k-i}u, \hspace{2mm} k\geq 2,
\end{equation*}
and
\begin{equation}
(\frac{\partial}{\partial t}-\Delta)|\nabla^ku|^2=-2|\nabla^{k+1} u|^2+ \sum_{i=0}^{k-2} \nabla^i Rm * \nabla^{k-i}u*\nabla^ku, \hspace{2mm} k\geq 2,
\end{equation}
where, as usual, for tensors $A$ and $B$, $A*B$ denotes a linear combination of contractions of the tensor product $A\otimes B$.

\vspace*{0.4cm}

\noindent{\bf Proof of Theorem 1.1}.

Let $G_1=(A_1a^2+u^2)|\nabla u|^2$, where $A_1$ is a positive constant to be chosen depending only on the dimension.
Using (2.1) we get
\begin{equation*}
(\frac{\partial}{\partial t}-\Delta)G_1=-2(A_1a^2+u^2)|\nabla^2 u|^2-2|\nabla u|^4+u\nabla u * \nabla u * \nabla^2 u.
\end{equation*}

On $PB_r(x_0, T)$, using our assumption  we have
\begin{equation*}
|u\nabla u * \nabla u * \nabla^2 u|\leq Ca|\nabla u|^2|\nabla^2 u|,
\end{equation*}
 where $C$ is a constant depending only on the dimension, so
\begin{equation*}
|u\nabla u * \nabla u * \nabla^2 u| \leq A_1a^2|\nabla^2 u|^2+|\nabla u|^4
\end{equation*}
for $A_1\geq \frac{1}{4}C^2$, and
\begin{equation*}
(\frac{\partial}{\partial t}-\Delta)G_1 \leq -|\nabla u|^4.
\end{equation*}
Choose $b_1=\frac{1}{(A_1+1)^2a^4}$, and let $F_1=b_1G_1$. Then
\begin{equation*}
(\frac{\partial}{\partial t}-\Delta)F_1 \leq -F_1^2.
\end{equation*}

Since $Ric\leq \frac{n-1}{r^2}$ on $PB_r(x_0, T)$, by Hamilton \cite{H13} (also see the proof of Corollary 5.2 in the Appendix), we can construct a  function $\Psi_1$ on     $\{(x,t)| x\in B(x_0, t, r), t\in [0,T]\}$       of the form
\begin{equation*}
\Psi_1=\frac{\alpha_1  r^2}{(r^2-d_t(x,x_0)^2)^2}
\end{equation*}
(where $\alpha_1$ is a positive constant depending only on the dimension) satisfying
\begin{equation}
 (\frac{\partial}{\partial t}-\Delta)\Psi_1 > -\Psi_1^2
\end{equation}
 on $\{(x,t)| x\in B(x_0, t, r), t\in (0,T]\}$ in the constructive comparison sense (for definition see the statement of Corollary 5.2).

  Now let
\begin{equation*}
   \Phi_1=\Psi_1+\frac{1}{t} =\frac{\alpha_1  r^2}{(r^2-d_t(x,x_0)^2)^2}+\frac{1}{t}
\end{equation*}
 on  $\{(x,t)| x\in B(x_0, t, r), t\in (0,T]\}$.
Then from (2.4) we immediately have
\begin{equation*}
(\frac{\partial}{\partial t}-\Delta)\Phi_1 > -\Phi_1^2
\end{equation*}
everywhere on $\{(x,t)| x\in B(x_0, t, r), t\in (0,T]\}$ in the constructive comparison sense.

Note that  $\Phi_1 \rightarrow \infty$ as $(x,t)$ tends to the parabolic boundary of $PB_r(x_0,T)$, but $F_1$ is bounded  on $PB_r(x_0,T)$  as $F_1$ is smooth on $M\times [0,T]$ and $PB_r(x_0,T)$ is compact.  So near the parabolic boundary of $PB_r(x_0, T)$ we have $F_1 < \Phi_1$.
 Using  a maximum principle argument as in the proof of Theorem 5.3 in the Appendix  we get that $F_1 < \Phi_1$  everywhere on $\{(x,t)| x\in B(x_0, t, r), t\in (0,T]\}$, and in particular,
\begin{equation*}
b_1A_1a^2|\nabla u|^2 \leq \frac{\alpha_1  r^2}{(r^2-d_t(x,x_0)^2)^2}+\frac{1}{t}.
\end{equation*}
On $PB_{ \frac{r}{2}}(x_0, T)\setminus \{(x, 0) | x \in M\}$ we have $d_t(x,x_0)\leq \frac{r}{2}$ and $r^2-d_t(x,x_0)^2\geq \frac{3}{4}r^2$, and the result follows by our choice of $b_1$.
\hfill{$\Box$}

\vspace*{0.4cm}

\noindent{\bf Remark}. Note that the comparison function $\Phi_1$ (in the proof above) blows up at the parabolic boundary of $PB_r(x_0, T)$, moreover it satisfies an inequality which is  opposite to the one satisfied by $F_1$, so $\Phi_1$ serves as a barrier for $F_1$, the latter being bounded on $PB_r(x_0, T)$.  Here we need the assumptions that the solution $u$ is smooth on $PB_r(x_0,T)$ and that $PB_r(x_0,T)$ is compact. But we do not need the completeness of the metrics $g(t)$ on $M$. Also note that in the statement of Theorem 1.1, if we assume in addition $|\nabla u| \leq \frac{a}{r}$ at $t=0$ in $\overline{B(x_0, 0, r)}$, then  we have $|\nabla u| \leq C_1\frac{a}{r}$ on $PB_\frac{r}{2}(x_0, T)$, because in this case we can choose $\Psi_1$ instead of $\Phi_1$ as  the (space-time) comparison function.

\vspace*{0.4cm}

\noindent{\bf Proof of Theorem 1.2}.

 On $PB_r(x_0, T)$, using (2.2) and our assumption on $|Rm|$, we have
\begin{equation*}
(\frac{\partial}{\partial t}-\Delta)|\nabla^2 u|^2 \leq -2|\nabla^3 u|^2+\frac{C}{r^2}|\nabla^2 u|^2,
\end{equation*}
where $C$ depends only on the dimension.

Let
\begin{equation*}
G_2=(A_2a^2(\frac{1}{r^2}+\frac{1}{t})+|\nabla u|^2)|\nabla^2 u|^2,
\end{equation*}
 where $A_2$ is a positive constant to be chosen depending only on the dimension. We have
\begin{equation*}
\begin{split}
&(\frac{\partial}{\partial t}-\Delta)G_2 \\
\leq &-\frac{A_2a^2}{t^2}|\nabla^2 u|^2-2|\nabla^2 u|^4+(A_2a^2(\frac{1}{r^2}+\frac{1}{t})+|\nabla u|^2)(-2|\nabla^3 u|^2+\frac{C}{r^2}|\nabla^2 u|^2)\\
&+C|\nabla u| |\nabla^2 u|^2 |\nabla^3 u|.
\end{split}
\end{equation*}
 On $PB_\frac{r}{2}(x_0,  T)\setminus \{(x, 0) | x \in M\}$ we have
\begin{equation*}
|\nabla u| \leq C_1a(\frac{1}{r}+\frac{1}{\sqrt{t}})
\end{equation*}
by Theorem 1.1, so

\begin{equation*}
(A_2a^2(\frac{1}{r^2}+\frac{1}{t})+|\nabla u|^2)\frac{C}{r^2}|\nabla^2 u|^2 \leq  \frac{1}{2}|\nabla^2 u|^4+C^2(A_2+2C_1^2)^2a^4\frac{1}{r^4}(\frac{1}{r^2}+\frac{1}{t})^2,
\end{equation*}
and
\begin{equation*}
C|\nabla u| |\nabla^2 u|^2 |\nabla^3 u|\leq \frac{1}{2}|\nabla^2 u|^4+A_2a^2(\frac{1}{r^2}+\frac{1}{t})|\nabla^3 u|^2
\end{equation*}
by choosing  $A_2\geq C^2C_1^2$.

\noindent Then
\begin{equation*}
\begin{split}
&(\frac{\partial}{\partial t}-\Delta)G_2 \\
\leq &-|\nabla^2 u|^4+C^2(A_2+2C_1^2)^2a^4\frac{1}{r^4}(\frac{1}{r^2}+\frac{1}{t})^2\\
\leq &-\frac{G_2^2}{(A_2+2C_1^2)^2a^4(\frac{1}{r^2}+\frac{1}{t})^2} +C^2(A_2+2C_1^2)^2a^4\frac{1}{r^4}(\frac{1}{r^2}+\frac{1}{t})^2.
\end{split}
\end{equation*}

Let $v=\frac{1}{r^2}+\frac{1}{t}$ and $F_2=\frac{b_2G_2}{v}$, where $b_2$ is a positive constant to  be chosen later.
On $PB_\frac{r}{2}(x_0,  T)\setminus \{(x, 0) | x \in M\}$ we have
\begin{equation*}
\begin{split}
&(\frac{\partial}{\partial t}-\Delta)F_2 \\
\leq &-\frac{F_2^2}{b_2(A_2+2C_1^2)^2a^4v}+b_2C^2(A_2+2C_1^2)^2a^4\frac{v}{r^4}+F_2\frac{1}{vt^2}\\
\leq &-\frac{F_2^2}{b_2(A_2+2C_1^2)^2a^4v}+b_2C^2(A_2+2C_1^2)^2a^4v^3+F_2v\\
\leq &-\frac{F_2^2}{2b_2(A_2+2C_1^2)^2a^4v}+b_2(A_2+2C_1^2)^2a^4(\frac{1}{2}+C^2)v^3,
\end{split}
\end{equation*}
 where in the last inequality we use
\begin{equation*}
F_2v \leq \frac{F_2^2}{2b_2(A_2+2C_1^2)^2a^4v}+\frac{b_2}{2}(A_2+2C_1^2)^2a^4v^3.
\end{equation*}
Choose  $b_2 = ((A_2+2C_1^2)^2a^4(2+C^2))^{-1}$.  Then we have
\begin{equation*}
(\frac{\partial}{\partial t}-\Delta)F_2  \leq -\frac{F_2^2}{v}+v^3.
\end{equation*}

Write $s=s(x,t)=d_t(x, x_0)$, and let
\begin{equation*}
\Psi_2(s)=\frac{\alpha_2r^2}{(\frac{r^2}{4}-s^2)^2}
\end{equation*}
on      $\{(x,t)| x\in B(x_0, t, \frac{r}{2}), t\in [0, T]\}$. Again as in  \cite{H13} (see also the proof of Corollary 5.2 below) we can choose positive constant $\alpha_2$ depending only on the dimension such that
\begin{equation*}
(\frac{\partial}{\partial t}-\Delta)\Psi_2 > -\Psi_2^2
 \end{equation*}
 everywhere on  $\{(x,t)| x\in B(x_0, t, \frac{r}{2}), t\in (0, T]\}$ in the constructive comparison sense.
Let
\begin{equation*}
\Phi_2=\beta \Psi_2^2+ \gamma \frac{1}{t^2}=\beta \frac{\alpha_2^2r^4}{(\frac{r^2}{4}-s^2)^4}+ \gamma \frac{1}{t^2}
\end{equation*}
on $\{(x,t)| x\in B(x_0, t, \frac{r}{2}), t\in (0, T]\}$, where $\beta$ and $\gamma$ are positive constants  to be chosen later.
We have
\begin{equation*}
\begin{split}
&(\frac{\partial}{\partial t}-\Delta)\Phi_2 \\
=&2\beta \Psi_2(\frac{\partial}{\partial t}-\Delta)\Psi_2-\frac{2\gamma}{t^3}-2\beta |\nabla \Psi_2|^2\\
> &-2\beta \Psi_2^3 -\frac{2\gamma}{t^3}-2\beta \Psi_2'(s)^2.
\end{split}
\end{equation*}
We will choose constants $\beta$ and $\gamma$  such that
\begin{equation*}
-2\beta \Psi_2^3 -\frac{2\gamma}{t^3}-2\beta \Psi_2'(s)^2 \geq -\frac{\Phi_2^2}{v}+v^3,
\end{equation*}
that is,
\begin{equation*}
\Phi_2^2 \geq 2\beta (\Psi_2^3 +  \Psi_2'(s)^2)v + \frac{2\gamma v}{t^3}+v^4.
\end{equation*}
On $\{(x,t)| x\in B(x_0, t, \frac{r}{2}), t\in (0, T]\}$  we have $s^2 < \frac{r^2}{4}$, and
\begin{equation*}
\Psi_2'(s)^2=\frac{16\alpha_2^2r^4s^2}{(\frac{r^2}{4}-s^2)^6}< \frac{4\alpha_2^2r^6}{(\frac{r^2}{4}-s^2)^6},
\end{equation*}
so it suffices to have
\begin{equation*}
\beta^2 \alpha_2^4 \frac{r^8}{(\frac{r^2}{4}-s^2)^8}+\frac{\gamma^2}{t^4} \geq 2\beta(\alpha_2^3+4\alpha_2^2) \frac{r^6}{(\frac{r^2}{4}-s^2)^6}v +\frac{2\gamma}{t^3}v+v^4.
\end{equation*}
Note that
\begin{equation*}
v^4\leq 8(\frac{1}{r^8}+\frac{1}{t^4}).
\end{equation*}

 Using the Young's  inequality
  \begin{equation*}
  y^3z\leq \frac{3}{4}(y^3)^\frac{4}{3}+\frac{1}{4}z^4=\frac{3}{4} y^4+\frac{1}{4}z^4
\end{equation*}
for $y, z \in \mathbb{R}$,  we get
 \begin{equation*}
  2\beta(\alpha_2^3+4\alpha_2^2)\frac{r^6}{(\frac{r^2}{4}-s^2)^6} \frac{1}{t}
  \leq \frac{3}{4}(2\beta (\alpha_2^3+4\alpha_2^2))^{\frac{4}{3}}\frac{r^8}{(\frac{r^2}{4}-s^2)^8}+\frac{1}{4t^4},
  \end{equation*}
 and
\begin{equation*}
\frac{1}{t^3} \frac{1}{r^2}  \leq \frac{3}{4t^4}+\frac{1}{4r^8}.
\end{equation*}
We also have
\begin{equation*}
\frac{r^6}{(\frac{r^2}{4}-s^2)^6}\frac{1}{r^2}\leq \frac{r^8}{16(\frac{r^2}{4}-s^2)^8}
\end{equation*}
and
\begin{equation*}
\frac{1}{r^8}\leq \frac{r^8}{4^8(\frac{r^2}{4}-s^2)^8}.
\end{equation*}

First choose $\gamma >0 $  such that
\begin{equation*}
\gamma^2 \geq \frac{7\gamma}{2}+\frac{33}{4}.
\end{equation*}
  Then choose $\beta >0$ depending only on the dimension  such that
\begin{equation*}
  \beta^2 \alpha_2^4 \geq \frac{3}{4}(2\beta (\alpha_2^3+4\alpha_2^2))^{\frac{4}{3}}+\frac{1}{8}\beta (\alpha_2^3+4\alpha_2^2)+\frac{1}{4^8}(\frac{\gamma}{2}+8).
\end{equation*}
With $\beta$ and $\gamma$ chosen this way  we have
\begin{equation*}
(\frac{\partial}{\partial t}-\Delta)\Phi_2 > -\frac{\Phi_2^2}{v}+v^3
\end{equation*}
 everywhere on  $\{(x,t)| x\in B(x_0, t, \frac{r}{2}), t\in (0, T]\}$ in the constructive comparison sense.

 Note that near the parabolic boundary of $PB_{\frac{r}{2}}(x_0, T)$ we have $F_2 < \Phi_2$.

\noindent Using a maximum principle argument as in the proof of Theorem 5.3 in the Appendix we get that $F_2 < \Phi_2$ everywhere on  $\{(x,t)| x\in B(x_0, t, \frac{r}{2}), t\in (0, T]\}$,  and in particular,
\begin{equation*}
b_2A_2a^2|\nabla^2 u|^2 \leq \beta \frac{\alpha_2^2r^4}{(\frac{r^2}{4}-s^2)^4}+ \gamma \frac{1}{t^2}.
\end{equation*}
On $PB_{ \frac{r}{4}}(x_0, T)\setminus \{(x, 0) | x \in M\}$ we have $s\leq \frac{r}{4}$ and $\frac{r^2}{4}-s^2 \geq \frac{3}{16}r^2$, and the result follows by our choice of $b_2$.
\hfill{$\Box$}

\vspace*{0.4cm}

\noindent{\bf Remark}. In the statement of Theorem 1.2, if we assume in addition $|\nabla u| \leq \frac{a}{r}$ and $|\nabla^2 u| \leq \frac{a}{r^2}$  at $t=0$ in $\overline{B(x_0, 0, r)}$, then  we have $|\nabla^2 u| \leq C_2\frac{a}{r^2}$ on $PB_\frac{r}{4}(x_0, T)$, because in this case we can choose $(A_2a^2\frac{1}{r^2}+|\nabla u|^2)|\nabla^2 u|^2$  as $G_2$ above, and $\beta \Psi_2^2$ instead of $\Phi_2$ as  the (space-time) comparison function.

\vspace*{0.4cm}

Similarly to Theorems 1.1 and 1.2 we have
\begin{thm} \label{thm 2.1}  Let $M$ be  a manifold  (without boundary) of dimension $n$.  Suppose $g(t)$ is a solution (not necessarily complete) to the Ricci flow  on $M \times [0, T]$ for some $T >0$.  Fix $x_0 \in M$ and $r>0$. Assume that  the parabolic cylinder $PB_r(x_0,T)$ is compact, and $|Rm|\leq \frac{1}{r^2}$ on  $PB_r(x_0,  T)$. Let $ u $ be a  smooth solution to the heat equation $(\frac{\partial}{\partial t}-\Delta_{g(t)})u=0$ coupled to the Ricci flow on $M \times [0, T]$. Suppose $|u|\leq a$ on $PB_r(x_0, T)$,  where $a$ is a positive constant. Then for any $k\geq 2$,
\begin{equation*}
|\nabla^k u| \leq C_ka(\frac{1}{r^{k}}+\frac{1}{t^{k/2}})    \hspace{2mm}  \text{on}   \hspace{2mm}    PB_\frac{r}{2^k}(x_0,  T)\setminus \{(x, 0) | x \in M\},
\end{equation*}
where the constant $C_k$ depends only on $k$ and the dimension.
\end{thm}

\noindent{\bf Proof}.
The proof is by induction.  On $PB_\frac{r}{2}(x_0,  T)\setminus \{(x, 0) | x \in M\}$ we have
\begin{equation*}
|\nabla u| \leq C_1a(\frac{1}{r}+\frac{1}{\sqrt{t}})
\end{equation*}
by Theorem 1.1. For $k=2$, the result  is exactly Theorem 1.2. Suppose on $PB_\frac{r}{2^i}(x_0,  T)\setminus \{(x, 0) | x \in M\}$ ($2\leq i \leq k$) we have
\begin{equation*}
|\nabla^i u| \leq C_ia(\frac{1}{r^i}+\frac{1}{t^{i/2}}),
\end{equation*}
where $C_i$ depends only on $i$ and the dimension. Let
\begin{equation*}
G_{k+1}=(A_{k+1}a^2(\frac{1}{r^{2k}}+\frac{1}{t^k})+|\nabla^k u|^2)|\nabla^{k+1} u|^2,
\end{equation*}
 where $A_{k+1} >1$ is a  constant to be chosen depending only on $k$ and the dimension.

Using (2.3) we have
\begin{equation*}
\begin{split}
&(\frac{\partial}{\partial t}-\Delta)G_{k+1} \\
=&(-\frac{A_{k+1}a^2k}{t^{k+1}}-2 |\nabla^{k+1}u|^2 +   \sum_{i=0}^{k-2} \nabla^i Rm * \nabla^{k-i}u*\nabla^ku) |\nabla^{k+1}u|^2\\
&+(A_{k+1}a^2(\frac{1}{r^{2k}}+\frac{1}{t^k})+|\nabla^k u|^2)(-2 |\nabla^{k+2}u|^2 + \sum_{i=0}^{k-1} \nabla^i Rm * \nabla^{k+1-i}u*\nabla^{k+1}u)\\
&+\nabla^k u * \nabla^{k+1} u * \nabla^{k+1} u  * \nabla^{k+2} u.
\end{split}
\end{equation*}

 Since $|Rm|\leq \frac{1}{r^2}$ on    $PB_r(x_0,  T)$ by assumption,       we have
\begin{equation*}
|\nabla^i Rm| \leq C_i'\frac{1}{r^2}(\frac{1}{r^i}+\frac{1}{t^{i/2}})
\end{equation*}
on  $PB_\frac{r}{2^i}(x_0,  T)\setminus \{(x, 0) | x \in M\}$ ($1\leq i \leq k-1$) by  Shi's local derivative estimates (see Theorem 5.3 below), where $C_i'$ depends only on $i$ and the dimension.

Below we will use $C$  to denote various constants depending only on $k$ and the dimension, which may be different from line to line.
On $PB_\frac{r}{2^k}(x_0,  T)\setminus \{(x, 0) | x \in M\}$,  we have
\begin{equation*}
|\sum_{i=0}^{k-2} \nabla^i Rm * \nabla^{k-i}u*\nabla^ku| \leq C\frac{a^2}{r^2}(\frac{1}{r^{2k}}+\frac{1}{t^k}),
\end{equation*}
\begin{equation*}
|\sum_{i=0}^{k-2} \nabla^i Rm * \nabla^{k-i}u*\nabla^ku| |\nabla^{k+1}u|^2 \leq   \frac{1}{3}|\nabla^{k+1}u|^4 +  C\frac{a^4}{r^4}(\frac{1}{r^{2k}}+\frac{1}{t^k})^2,
\end{equation*}
  \begin{equation*}
 |\sum_{i=0}^{k-1} \nabla^i Rm * \nabla^{k+1-i}u*\nabla^{k+1}u| \leq C(\frac{1}{r^2}|\nabla^{k+1}u|^2+\frac{a^2}{r^2}(\frac{1}{r^{2(k+1)}}+\frac{1}{t^{k+1}})),
 \end{equation*}
 \begin{equation*}
\begin{split}
 &|(A_{k+1}a^2(\frac{1}{r^{2k}}+\frac{1}{t^k})+|\nabla^k u|^2)  \sum_{i=0}^{k-1} \nabla^i Rm * \nabla^{k+1-i}u*\nabla^{k+1}u| \\
\leq &\frac{1}{3} |\nabla^{k+1}u|^4 + C (A_{k+1}+2C_k^2)^2\frac{a^4}{r^4}(\frac{1}{r^{2k}}+\frac{1}{t^k})^2 \\
&+ C(A_{k+1}+2C_k^2)\frac{a^4}{r^2}(\frac{1}{r^{2k}}+\frac{1}{t^k})(\frac{1}{r^{2(k+1)}}+\frac{1}{t^{k+1}}),
\end{split}
\end{equation*}
and
\begin{equation*}
 |\nabla^k u * \nabla^{k+1} u * \nabla^{k+1} u  * \nabla^{k+2} u| \leq \frac{1}{3}|\nabla^{k+1} u|^4+ A_{k+1}a^2(\frac{1}{r^{2k}}+\frac{1}{t^k})|\nabla^{k+2}u|^2
\end{equation*}
 by choosing $A_{k+1}$  sufficiently large (compared to $C_k^2$).

\noindent So
\begin{equation*}
\begin{split}
&(\frac{\partial}{\partial t}-\Delta)G_{k+1} \\
\leq &- |\nabla^{k+1}u|^4 +  C(A_{k+1}+2C_k^2)^2\frac{a^4}{r^4}(\frac{1}{r^{2k}}+\frac{1}{t^k})^2 \\
&+C(A_{k+1}+2C_k^2)\frac{a^4}{r^2}(\frac{1}{r^{2k}}+\frac{1}{t^k})(\frac{1}{r^{2(k+1)}}+\frac{1}{t^{k+1}}) \\
\leq & - |\nabla^{k+1}u|^4 +  C (A_{k+1}+2C_k^2)^2\frac{a^4}{r^4}(\frac{1}{r^{4k}}+\frac{1}{t^{2k}}) \\
&+C(A_{k+1}+2C_k^2)\frac{a^4}{r^2}(\frac{1}{r^{2(2k+1)}}+\frac{1}{t^{2k+1}})\\
\leq &- |\nabla^{k+1}u|^4+ C(A_{k+1}+2C_k^2)^2\frac{a^4}{r^2}(\frac{1}{r^{2(2k+1)}}+\frac{1}{t^{2k+1}})\\
\leq &- \frac{G_{k+1}^2}{(A_{k+1}+2C_k^2)^2a^4(\frac{1}{r^{2k}}+\frac{1}{t^k})^2}+ C(A_{k+1}+2C_k^2)^2\frac{a^4}{r^2}(\frac{1}{r^{2(2k+1)}}+\frac{1}{t^{2k+1}}).
\end{split}
\end{equation*}

Let $v=\frac{1}{r^2}+\frac{1}{t}$ and $F_{k+1}=\frac{b_{k+1}G_{k+1}}{v^k}$ (cf. the first line on p.198 in \cite{CZ}),  where $b_{k+1}$ is a positive constant to  be chosen later. Then on $PB_\frac{r}{2^k}(x_0,  T)\setminus \{(x, 0) | x \in M\}$ we have
\begin{equation*}
\begin{split}
&(\frac{\partial}{\partial t}-\Delta)F_{k+1} \\
\leq &- \frac{F_{k+1}^2}{b_{k+1}(A_{k+1}+2C_k^2)^2a^4v^k}+ b_{k+1}C(A_{k+1}+2C_k^2)^2\frac{a^4}{r^2}v^{k+1}+kF_{k+1}v\\
\leq &- \frac{F_{k+1}^2}{2b_{k+1}(A_{k+1}+2C_k^2)^2a^4v^k}+ b_{k+1}(C+2k^2)(A_{k+1}+2C_k^2)^2a^4v^{k+2}.
\end{split}
\end{equation*}
Choosing $b_{k+1}=\frac{1}{(C+2k^2)(A_{k+1}+2C_k^2)^2a^4}$, we get
\begin{equation*}
(\frac{\partial}{\partial t}-\Delta)F_{k+1}\leq -\frac{1}{v^k}F_{k+1}^2+v^{k+2}.
\end{equation*}

Write $s=s(x,t)=d_t(x, x_0)$, and let
\begin{equation*}
\Psi_{k+1}(s)=\frac{\alpha_{k+1}r^2}{(\frac{r^2}{4^k}-s^2)^2}
\end{equation*}
on      $\{(x,t)| x\in B(x_0, t, \frac{r}{2^k}), t\in [0, T]\}$. Again as in  \cite{H13} (see also the proof of Corollary 5.2 below) we can choose constant $\alpha_{k+1}>0$ depending only on the dimension and $k$ such that
\begin{equation*}
(\frac{\partial}{\partial t}-\Delta)\Psi_{k+1} > -\Psi_{k+1}^2
 \end{equation*}
 everywhere on  $\{(x,t)| x\in B(x_0, t, \frac{r}{2^k}), t\in (0, T]\}$ in the constructive comparison sense.
Let
\begin{equation*}
\Phi_{k+1}=\beta_{k+1} \Psi_{k+1}^{k+1}+ \gamma_{k+1} \frac{1}{t^{k+1}}=\beta_{k+1} \frac{\alpha_{k+1}^{k+1}r^{2(k+1)}}{(\frac{r^2}{4^k}-s^2)^{2(k+1)}}+ \gamma_{k+1} \frac{1}{t^{k+1}}
\end{equation*}
on $\{(x,t)| x\in B(x_0, t, \frac{r}{2^k}), t\in (0, T]\}$, where $\beta_{k+1}$ and $\gamma_{k+1}$ are positive constants  to be chosen later.

We have
\begin{equation*}
\begin{split}
 &(\frac{\partial}{\partial t}-\Delta)\Phi_{k+1} \\
= & \beta_{k+1}(k+1)\Psi_{k+1}^k(\frac{\partial}{\partial t}-\Delta)\Psi_{k+1}-(k+1)\frac{\gamma_{k+1}}{t^{k+2}}-k(k+1)\beta_{k+1}\Psi_{k+1}^{k-1}|\nabla \Psi_{k+1}|^2\\
> & -\beta_{k+1}(k+1)\Psi_{k+1}^{k+2}-(k+1)\frac{\gamma_{k+1}}{t^{k+2}}-k(k+1)\beta_{k+1}\Psi_{k+1}^{k-1}\Psi_{k+1}'(s)^2.
\end{split}
\end{equation*}
 We will choose constants $\beta_{k+1}$ and $\gamma_{k+1}$  such that
\begin{equation*}
\begin{split}
& -\beta_{k+1}(k+1)\Psi_{k+1}^{k+2}-(k+1)\frac{\gamma_{k+1}}{t^{k+2}}-k(k+1)\beta_{k+1}\Psi_{k+1}^{k-1}\Psi_{k+1}'(s)^2 \\
 \geq & -\frac{\Phi_{k+1}^2}{v^k}+v^{k+2},
\end{split}
\end{equation*}
that is,
\begin{equation*}
\Phi_{k+1}^2 \geq \beta_{k+1}(k+1)(\Psi_{k+1}^{k+2} + k\Psi_{k+1}^{k-1} \Psi_{k+1}'(s)^2)v^k + (k+1)\frac{\gamma_{k+1} v^k}{t^{k+2}}+v^{2(k+1)}.
\end{equation*}
On $\{(x,t)| x\in B(x_0, t, \frac{r}{2^k}), t\in (0, T]\}$  we have $s^2 < \frac{r^2}{4^k}$, and
\begin{equation*}
\Psi_{k+1}'(s)^2=\frac{16\alpha_{k+1}^2r^4s^2}{(\frac{r^2}{4^k}-s^2)^6}< \frac{\alpha_{k+1}^2r^6}{4^{k-2}(\frac{r^2}{4^k}-s^2)^6},
\end{equation*}
so it suffices to have
\begin{equation*}
\begin{split}
& \beta_{k+1}^2 \alpha_{k+1}^{2(k+1)} \frac{r^{4(k+1)}}{(\frac{r^2}{4^k}-s^2)^{4(k+1)}}+\frac{\gamma_{k+1}^2}{t^{2(k+1)}} \\
\geq & \beta_{k+1}(k+1)(\alpha_{k+1}^{k+2}+\frac{k}{4^{k-2}}\alpha_{k+1}^{k+1}) \frac{r^{2(k+2)}}{(\frac{r^2}{4^k}-s^2)^{2(k+2)}}v^k +(k+1)\frac{\gamma_{k+1}}{t^{k+2}}v^k+v^{2(k+1)}.
\end{split}
\end{equation*}
Note that
\begin{equation*}
v^k\leq 2^{k-1}(\frac{1}{r^{2k}}+\frac{1}{t^k})
\end{equation*}
and
\begin{equation*}
v^{2(k+1)}\leq 2^{2k+1}(\frac{1}{r^{4(k+1)}}+\frac{1}{t^{2(k+1)}}).
\end{equation*}

 Using the Young's inequality
\begin{equation*}
\begin{split}
      y^{k+2}z^k\leq & \frac{k+2}{2(k+1)}(y^{k+2})^\frac{2(k+1)}{k+2}+\frac{k}{2(k+1)}(z^k)^{\frac{2(k+1)}{k}} \\
    =& \frac{k+2}{2(k+1)}y^{2(k+1)}+\frac{k}{2(k+1)}z^{2(k+1)}
\end{split}
\end{equation*}
 for $y, z \in \mathbb{R}$,  we get
\begin{equation*}
\begin{split}
  & 2^{k-1}\beta_{k+1}(k+1)(\alpha_{k+1}^{k+2}+\frac{k}{4^{k-2}}\alpha_{k+1}^{k+1}) \frac{r^{2(k+2)}}{(\frac{r^2}{4^k}-s^2)^{2(k+2)}}\frac{1}{t^k} \\
\leq & \frac{k+2}{2(k+1)}[2^{k-1}\beta_{k+1}(k+1)(\alpha_{k+1}^{k+2}+\frac{k}{4^{k-2}}\alpha_{k+1}^{k+1})]^\frac{2(k+1)}{k+2}  \frac{r^{4(k+1)}}{(\frac{r^2}{4^k}-s^2)^{4(k+1)}} \\
&+\frac{k}{2(k+1)}\frac{1}{t^{2(k+1)}},
\end{split}
\end{equation*}
and
\begin{equation*}
\frac{1}{t^{k+2}} \frac{1}{r^{2k}} \leq \frac{k+2}{2(k+1)}\frac{1}{t^{2(k+1)}}+\frac{k}{2(k+1)}\frac{1}{r^{4(k+1)}}.
\end{equation*}
We also have
\begin{equation*}
\frac{r^{2(k+2)}}{(\frac{r^2}{4^k}-s^2)^{2(k+2)}}\frac{1}{r^{2k}}\leq \frac{r^{4(k+1)}}{4^{2k^2}(\frac{r^2}{4^k}-s^2)^{4(k+1)}}
\end{equation*}
and
\begin{equation*}
\frac{1}{r^{4(k+1)}}\leq \frac{r^{4(k+1)}}{4^{4k(k+1)}(\frac{r^2}{4^k}-s^2)^{4(k+1)}}.
\end{equation*}

First choose $\gamma_{k+1} >0 $  depending only on $k$ such that
\begin{equation*}
\gamma_{k+1}^2 \geq 2^{k-1}(k+1)(1+\frac{k+2}{2(k+1)})\gamma_{k+1}+2^{2k+1}+\frac{k}{2(k+1)}.
\end{equation*}
  Then choose $\beta_{k+1} >0$ depending only on the dimension and $k$ such that
\begin{equation*}
\begin{split}
  \beta_{k+1}^2 \alpha_{k+1}^{2(k+1)}  &\geq   \frac{k+2}{2(k+1)}[2^{k-1}\beta_{k+1}(k+1)(\alpha_{k+1}^{k+2}+\frac{k}{4^{k-2}}\alpha_{k+1}^{k+1})]^\frac{2(k+1)}{k+2} \\
  &+ \beta_{k+1}\frac{k+1}{2^{4k^2-k+1}} (\alpha_{k+1}^{k+2}+\frac{k}{4^{k-2}}\alpha_{k+1}^{k+1}) +  \frac{1}{2^{8k^2+7k+2}}(k\gamma_{k+1}+2^{k+3}).
  \end{split}
\end{equation*}
With $\beta_{k+1}$  and $\gamma_{k+1}$ chosen this way we have
\begin{equation*}
(\frac{\partial}{\partial t}-\Delta)\Phi_{k+1} > -\frac{\Phi_{k+1}^2}{v^k}+v^{k+2}
\end{equation*}
 everywhere on  $\{(x,t)| x\in B(x_0, t, \frac{r}{2^k}), t\in (0, T]\}$ in the constructive comparison sense.

 Note that near the parabolic boundary of $PB_{\frac{r}{2^k}}(x_0, T)$ we have
 $F_{k+1} < \Phi_{k+1}$.
Using a maximum principle argument as in the proof of Theorem 5.3 below we get that $F_{k+1} < \Phi_{k+1}$ everywhere on
$\{(x,t)| x\in B(x_0, t, \frac{r}{2^k}), t\in (0, T]\}$,  and in particular,
\begin{equation*}
\frac{1}{2^{k-1}}b_{k+1}A_{k+1}a^2|\nabla^{k+1} u|^2 \leq \beta_{k+1} \frac{\alpha_{k+1}^{k+1}r^{2(k+1)}}{(\frac{r^2}{4^k}-s^2)^{2(k+1)}}+ \gamma_{k+1} \frac{1}{t^{k+1}}.
\end{equation*}
On $PB_{ \frac{r}{2^{k+1}}}(x_0, T)\setminus \{(x, 0) | x \in M\}$ we have $s\leq \frac{r}{2^{k+1}}$ and $\frac{r^2}{4^k}-s^2 \geq \frac{3}{4^{k+1}}r^2$, and the result follows by our choice of $b_{k+1}$.
\hfill{$\Box$}

\vspace*{0.4cm}

\noindent{\bf Remark}.  In the statement of Theorem 2.1, when $k\geq 3$, if we assume in addition   $|\nabla^i Rm| \leq \frac{1}{r^{2+i}}$    for $1\leq i \leq k-2$ and       $|\nabla^i u| \leq \frac{a}{r^i}$ for $1 \leq i \leq k$   at $t=0$ in $\overline{B(x_0, 0, r)}$, then  we have $|\nabla^i u| \leq C_i\frac{a}{r^i}$ on $PB_\frac{r}{2^i}(x_0, T)$ for $1 \leq i \leq k$, where $C_i$ is a constant depending only on $i$ and the dimension. The proof can be adapted from that of Theorem 2.1 by using Lu's modified  Shi estimates (see Theorem 5.4 below),  choosing $(A_ia^2\frac{1}{r^{2(i-1)}}+|\nabla^{i-1} u|^2)|\nabla^{i} u|^2$ as $G_i$ above,  and choosing $\beta_{i} \Psi_{i}^{i}$ instead of $\Phi_{i}$ ($1\leq i \leq k$) as  the (space-time) comparison function.

\vspace*{0.4cm}

Of course we can state Theorem 2.1 for all $k\geq 1$ and use Theorem 1.1 instead of Theorem 1.2 as the beginning of the induction, so the proof of Theorem 1.2  can be omitted. But we prefer to reserve  it since we do not need to use Shi's local derivative estimates in the  proof of Theorem 1.2 (in contrast to the case $k\geq 3$), and moreover, it serves as a guide for the proof of Theorem 2.1.  Note also that in the conclusion of Theorem 2.1 we can replace  $PB_\frac{r}{2^k}(x_0,  T)$ by  $PB_\frac{r}{2}(x_0,  T)$, of course, then the constant $C_k$ will be different.

\section{ Some applications of Theorems 1.1 and  1.2}

Using Theorem 1.1 we can extend an estimate in Theorem 3.3 in Zhang \cite{Z06} and Theorem 5.1 in Cao-Hamilton \cite{CH} to a more general situation.
\begin{prop} \label{prop 3.1} (cf. Zhang \cite{Z06}, Cao-Hamilton \cite{CH})
Let $(M, (g(t))_{t\in (0,T)})$ be a complete solution to the Ricci flow with bounded Ricci curvature on any compact time subinterval. Let $0< u\leq a$ be a solution to the heat equation
$\frac{\partial}{\partial t} u=\Delta_{g(t)}u$ coupled to the Ricci flow on $M\times (0,T)$, where $a$ is a positive constant.  Then
\begin{equation*}
\frac{|\nabla u(x,t)|}{u(x,t)}\leq \sqrt{\frac{1}{t}}\sqrt{\ln \frac{a}{u(x,t)}}   \hspace{2mm} on  \hspace{2mm}    M \times (0,T).
\end{equation*}
\end{prop}
\noindent {\bf Proof}.   Compare the proof of Lemma 6.3 in \cite{CTY}. Clearly we can assume that $T < \infty$; otherwise we only need to restrict to every finite time subinterval $(0, T')$. Then  using Theorem 1.1 and a standard trick (cf \cite{H93}, \cite{B}) we can reduce the proof in the general case to the case that $g(t)$ extends smoothly up to $t=0$ with  $\sup_{(x,t)\in M \times [0,T]} |Ric|< \infty$  and $\sup_{(x,t)\in M \times [0,T]} |\nabla  u|< \infty$. The reason is as follows: Fix $(x_0, t_0)\in M \times (0,T)$. Choose a small $\varepsilon >0$ such that $t_0 \in (\varepsilon, T-2\varepsilon)$. By assumption
\begin{equation*}
\sup_{(x,t)\in M \times [\frac{\varepsilon}{2},T-\varepsilon]} |Ric|< \infty
\end{equation*}
and $0<u\leq a$. By Theorem 1.1 we have
\begin{equation*}
\sup_{(x,t)\in M \times [\varepsilon,T-\varepsilon]} |\nabla  u|< \infty.
\end{equation*}
Let $\tilde{g}(t)=g(t+\varepsilon)$ and $\tilde{u}(t)= u(t+\varepsilon)$, $t\in [0,T-2\varepsilon]$. Note that $\tilde{g}(t)$ is also a solution to the Ricci flow, and $\tilde{u}(t)$ is also a solution to the heat equation coupled to the Ricci flow $\tilde{g}(t)$. Now  $\sup_{(x,t)\in M \times [0,T-2\varepsilon]} |\widetilde{Ric}|_{\tilde{g}(t)}< \infty$, $0< u \leq a$, and $\sup_{(x,t)\in M \times [0,T-2\varepsilon]} |\nabla_{\tilde{g}(t)}  \tilde{u}|_{\tilde{g}(t)}< \infty$.  Suppose in this case we have
\begin{equation*}
\frac{|\nabla_{\tilde{g}(t)} \tilde{u}(x,t)|_{\tilde{g}(t)}}{\tilde{u}(x,t)}\leq \sqrt{\frac{1}{t}}\sqrt{\ln \frac{a}{\tilde{u}(x,t)}},  \hspace{4mm}  t\in (0,T-2\varepsilon].
\end{equation*}
 In particular the above inequality holds at $(x_0, t_0)$. Now letting $\varepsilon \rightarrow 0$ we get the desired inequality.

Now let $v=u+\delta$, where $\delta$ is a positive constant. Then $\delta \leq v \leq a+\delta$ is a solution to the heat equation coupled to the Ricci flow with bounded gradient. Now as in \cite{Z06} and \cite{CH} we have

\begin{equation*}
(\frac{\partial}{\partial t}-\Delta_{g(t)})(t\frac{|\nabla v|^2}{v}-v\ln \frac{a+\delta}{v}) \leq 0, \hspace{4mm}   t\in [0,T].
\end{equation*}

Since  $\delta \leq v \leq a+\delta$ and $\sup_{M \times [0,T]}|\nabla v| < \infty$, we have $\sup_{(x,t)\in M \times [0,T]} t\frac{|\nabla v|^2}{v} < \infty$ and $\sup_{M \times [0,T]} v\ln \frac{a+\delta}{v} < \infty$. We also have
$\sup_{(x,t)\in M \times [0,T]} |Ric|< \infty$.  So with the help of Bishop-Gromov volume comparison theorem, by the maximum principle (Theorem 12.22 in \cite{C+08})  we have
\begin{equation*}
t\frac{|\nabla v|^2}{v}-v\ln \frac{a+\delta}{v}\leq 0
\end{equation*}
everywhere, since it is true when $t=0$.   Now letting $\delta \rightarrow 0$ and we are done.   \hfill{$\Box$}

\vspace*{0.4cm}

\noindent{\bf Remark}.  Theorem 3.3 in \cite{Z06} is stated for complete manifolds and does not impose any curvature bound,  but Theorem 6.5.1 in \cite{Z11} assumes the curvature is uniformly bounded. In both places   the details on justifying  the use of the maximum principle are not supplied. Moreover, our  statement  is slightly more general than that in Theorem 6.5.1 of \cite{Z11} in that we do not assume the Ricci flow is defined at $t=0$, so in our case  the curvature is not necessarily uniformly bounded.  Actually we only assume the Ricci curvature is bounded on any compact time subinterval.

\vspace*{0.4cm}

Using Theorems 1.1 and 1.2 we also extend Lemma 3.1 in Bamler-Zhang \cite{BZ} to the noncompact case.
\begin{prop} \label{prop 3.2} Let $(M, (g(t))_{t\in (0,T)})$ be a complete solution to the Ricci flow with bounded curvature on any compact time subinterval. Let $0< u\leq a$ be a solution to the heat equation $\frac{\partial}{\partial t} u=\Delta_{g(t)}u$ coupled to the Ricci flow  on $M\times (0,T)$, where $a$ is a positive constant. Then
\begin{equation*}
(|\Delta u|+\frac{|\nabla u|^2}{u}-aR)(x,t)\leq \frac{Ba}{t}     \hspace{2mm} on  \hspace{2mm}    M \times (0,T),
\end{equation*}
where the constant B depends only on the dimension.
\end{prop}
\noindent {\bf Proof}. As before  we can assume that $T < \infty$. Then as in the proof of Proposition 3.1 we can reduce the proof in the general case to the case that  $g(t)$ extends smoothly up to $t=0$ with $\sup_{(x,t)\in M \times [0,T]} | Rm|< \infty$ and  $\sup_{(x,t)\in M \times [0,T]} |\nabla ^k u|< \infty$, $k= 1, 2$: Fix $(x_0, t_0)\in M \times (0,T)$. Choose a small $\varepsilon >0$ such that $t_0 \in (\varepsilon, T-2\varepsilon)$. By Theorem 1.1 and  Theorem 1.2 we have   $\sup_{(x,t)\in M \times [\varepsilon,T-\varepsilon]} |\nabla ^k u|< \infty$ for $k=1, 2$. Let $\tilde{g}(t)=g(t+\varepsilon)$ and $\tilde{u}(t)= u(t+\varepsilon)$, $t\in [0,T-2\varepsilon]$. Then  $\sup_{(x,t)\in M \times [0,T-2\varepsilon]} |\widetilde{Rm}|_{\tilde{g}(t)}< \infty$, $\sup_{(x,t)\in M \times [0,T-2\varepsilon]} |\nabla_{\tilde{g}(t)} ^k \tilde{u}|_{\tilde{g}(t)}< \infty$ for $k=1, 2$.  Suppose in this case we have
\begin{equation*}
(|\Delta_{\tilde{g}(t)} \tilde{u}|+\frac{|\nabla_{\tilde{g}(t)} \tilde{u}|_{\tilde{g}(t)}^2}{\tilde{u}}-a\tilde{R})(x,t)\leq \frac{Ba}{t},  \hspace{4mm}  t\in (0,T-2\varepsilon],
\end{equation*}
where the constant $B$ depends only on the dimension. In particular the above inequality holds at $(x_0, t_0)$. Now letting $\varepsilon \rightarrow 0$ we get the desired inequality.

Also note that by the same trick of replacing $u$ by $u+\delta$ and letting $\delta \rightarrow 0$ as in the proof of Proposition 3.1 we can assume that $u\geq \delta >0$.

By rescaling we may assume that $a=1$. Let $L_1=-\Delta u+\frac{|\nabla u|^2}{u}-R$, and choose $B>0$ with $\frac{B+e^{-2}}{B^2}=\frac{1}{n}$, then as in the proof of Lemma 3.1 in \cite{BZ}, we have
\begin{equation*}
(\frac{\partial}{\partial t}-\Delta_{g(t)})(L_1-\frac{B}{t}) \leq -\frac{1}{n}(L_1+\frac{B}{t})(L_1-\frac{B}{t})
\end{equation*}
for $t \in (0,T]$.

Now given any $\varepsilon >0$, $C>0$, let $\varphi(x,t)=\varepsilon e^{At}f(x)$ be a positive function as in the proof of Lemma 5.2 in \cite{H93}, which satisfies $f(x)\rightarrow \infty$ as $x\rightarrow \infty$ and $(\frac{\partial}{\partial t}-\Delta_{g(t)})\varphi > C\varphi$. So we have
\begin{equation}
(\frac{\partial}{\partial t}-\Delta_{g(t)})(L_1-\frac{B}{t}-\varphi) \leq -\frac{1}{n}(L_1+\frac{B}{t})(L_1-\frac{B}{t})- C\varphi
\end{equation}
for $t \in (0,T]$. We claim $L_1-\frac{B}{t}-\varphi <0$ for $t \in (0,T]$. Note that this is true for $t>0$ sufficiently small by our assumption on $|Rm|$, $u$ and $|\nabla ^k u|$, $k=1,2$. Suppose it is not true for some large $t$. Then there exist the first time $t_0$ and a point $x_0$ such that $L_1(x_0,t_0)-\frac{B}{t_0}-\varphi(x_0,t_0) =0$ since $f(x)\rightarrow \infty$ as $x\rightarrow \infty$. Now  at $(x_0,t_0)$ we have  $\frac{\partial}{\partial t}(L-\frac{B}{t}-\varphi)\geq 0$, and $\Delta (L-\frac{B}{t}-\varphi) \leq 0$. This contradicts (3.1), since  at $(x_0,t_0)$ the RHS of (3.1) $<0$.  Now letting $\varepsilon \rightarrow 0$ we get $L_1\leq \frac{B}{t}$.

Let  $L_2=\Delta u+\frac{|\nabla u|^2}{u}-R $, and choose $B >0$ with $B^{-1}+\frac{1+\frac{4}{n}}{e^2}B^{-2}=\frac{1}{2n}$.  As in the proof of Lemma 3.1 in \cite{BZ} we have
\begin{equation*}
(\frac{\partial}{\partial t}-\Delta_{g(t)})(L_2-\frac{B}{t}) \leq -\frac{1}{2n}(L_2+\frac{B}{t})(L_2-\frac{B}{t})
\end{equation*}
for $t \in (0,T]$.  Arguing  as above we get the desired inequality for $L_2$.
\hfill{$\Box$}

\section{Perelman's W-entropy on noncompact manifolds}

Now as in for example \cite{Ku1}, \cite{Z12b}, \cite{RV} and \cite{L}, we consider Perelman's W-entropy  (see \cite{P})
\begin{equation*}
W(g,v, \tau)=\int_M [\tau(4|\nabla v|^2+Rv^2)-v^2 \ln v^2-\frac{n}{2}(\ln 4\pi \tau)v^2-nv^2]dg
\end{equation*}
on a complete noncompact Riemannian manifold $(M, g)$ of dimension $n$, where  $v\in W^{1,2}(M,g)$,  $\tau >0$ is a parameter, and $dg$ denotes the volume element of the metric $g$ as in \cite{Z12b}. Note that by Theorem 3.1 in \cite{He}, $ W^{1,2}(M,g)=W_0^{1,2}(M,g)$. Let
\begin{equation*}
\mu(g,\tau)=\inf \{W(g, v, \tau)  \hspace{2mm} | \hspace{2mm}  v\in C_0^\infty(M), ||v||_{L^2(M,g)}=1\}.
\end{equation*}
For $(M, g)$ with Ricci curvature  bounded below  and injectivity radius bounded away from 0, we have that $\mu(g,\tau)$ is  finite; for a proof see for example \cite{RV}.

The following  proposition is a  slight improvement of some results in \cite{Ku1}, \cite{Z12b} and \cite{L} which in turn extend the entropy formula in Perelman \cite{P} to the noncompact case. The improvement is partially on lowering the order of derivatives of the curvature tensor  which are required to be uniformly bounded to guarantee the equality (4.1) below, compare Corollary 4.1 in \cite{Z12b}, and for (4.1) we also  need not the condition on the injectivity radius  which is imposed in \cite{Z12b}. Moreover we allow the function $v_T$ below with slightly less constraints. Note that in the proof of Theorem 16 in \cite{Ku1} Kuang only considers  $v_T$ with compact support; see Remark 17 there.
The $v_T$ considered in Corollary 4.1 in \cite{Z12b} is also very special. We also clarify a key point in the proof of Theorem 16 in \cite{Ku1}.

\begin{prop} \label{prop 4.1} (cf. \cite{Ku1}, \cite{Z12b} and \cite{L})
Let  $(M, g_0)$ be a complete noncompact Riemannian manifold with bounded curvature such that  $\sup_M|\nabla Rm_{g_0}|< \infty$ and $\sup_M|\nabla^2 Rm_{g_0}|< \infty$.  Let $(M, (g(t))_{t \in [0,T]})$ be the complete solution to the Ricci flow  with $\sup_{M\times [0,T]}|Rm|< \infty$ and with $g(0)=g_0$.
Let $v_T \in C^\infty(M)$ with $ |v_T(x)| \leq Ae^{-ad_{T}(x,x_0)^2}$ for any $x\in M$ and $\int_M v_T^2dg(T)=1$, where $A$ and $a$ are  positive constants, and $x_0$ is a fixed point in $M$. Assume that $u$ is the solution to the conjugate heat equation coupled to the Ricci flow, $ \frac{\partial u}{\partial t}+\Delta_{g(t)} u-Ru=0$,  with $u(x,T)=v_T(x)^2$. Let $v(x,t)=\sqrt{u(x,t)}$ and $\tau (t)= T-t$ for $t\in [0,T)$.
Then  $W(g(t),v(\cdot,t), \tau(t))$ is finite and
\begin{equation}
\frac{d}{dt}W(g(t),v(\cdot,t), \tau(t))=2\tau (t) \int_M |Ric-Hess \ln u-\frac{1}{2\tau (t)}g(t)|^2udg(t)
\end{equation}
for $t\in [0,T)$. Consequently  if we assume in addition that the injectivity radius of  $g_0$ is  bounded away from 0,  we have
\begin{equation*}
\mu (g(t_1), \tau(t_1)) \leq \mu (g(t_2), \tau(t_2))
\end{equation*}
for $0\leq t_1 < t_2 < T$.
\end{prop}
\noindent{\bf Proof}.  In the proof below we need $\sup_{M\times [0,T]}|\nabla  Ric|< \infty$  and
 $\sup_{M\times [0,T]}|\Delta  R|< \infty$ in addition to $\sup_{M\times [0,T]}|Rm|< \infty$, which is implied by our assumptions  via Lu's modified version (see \cite{LT}, Theorem 3.29 in \cite{MT}, Theorem 14.16 in \cite{C+08}, \cite{H13}, and Theorem 5.4 below) of Shi's  derivative estimates.

Now let $v_T$ and $u$ be as in the statement of Proposition 4.1.  As in \cite{P} (see also \cite{Z12b}) let
\begin{equation*}
P(u)=\tau (-2\Delta u+\frac{|\nabla u|^2}{u}+Ru)-u\ln u-\frac{n}{2}(\ln 4\pi \tau)u-nu.
\end{equation*}
 By Proposition 9.1 in \cite{P} we have
 \begin{equation}
 H^*P(u)=2\tau (t) |Ric-Hess \ln u-\frac{1}{2\tau (t)}g(t)|^2u,
 \end{equation}
where $H^*=\frac{\partial}{\partial t}+\Delta-R$.

Below we will analyse each term in $P(u)$.
 As in \cite{Z12b} let $G(x,t;y,T)$  denote the fundamental solution of the conjugate heat equation coupled to the Ricci flow.
Using Corollary 5.6 in \cite{CTY} (see also Corollary 26.26 in \cite{C+10}) we have
\begin{equation*}
G(x,t;y,T)\leq \frac{\alpha}{|B(y,t,\sqrt{\frac{T-t}{2}})|_{g(t)}}e^{-\beta \frac{d_t(x,y)^2}{T-t}}
\end{equation*}
 for $t\in [0,T)$, and $\alpha$ and $\beta$ are positive constants  independent of $x$, $y$ and $t$.  On the other hand, by \cite{CLY} and  \cite{CGT} $|B(y,t,\sqrt{\frac{T-t}{2}})|_{g(t)}$
has at worst linear exponential decay as $y$ goes to infinity, that is,
\begin{equation*}
|B(y,t,\sqrt{\frac{T-t}{2}})|_{g(t)}\geq \gamma e^{-\delta d_t(y,x_0)},
\end{equation*}
 where $x_0$ is a fixed point in $M$, $\gamma$ and $\delta$ are positive constants independent of $y$, and $\delta$ is also independent of $x_0$. Now from the formula
\begin{equation*}
u(x,t)=\int_M G(x,t;y,T)u(y,T)dg(T)
\end{equation*}
 we can easily show that $u(\cdot,t)$ also has quadratic exponential decay for any $t\in [0,T)$ as $u(\cdot,T)$ does; compare Step 1 in the proof of Corollary 4.1 in \cite{Z12b}. It follows that  $|u\ln u|$ also has quadratic exponential decay.  Combining this with Theorem 10 in \cite{EKNT} we see that
$\tau \frac{|\nabla u|^2}{u}$  also has this decay.  Moreover from (3.27) in \cite{Ku1} we have $\int_M \Delta u dg(t)=0$.  Using Lemma 4.1 in \cite{CTY} and the decay property of $u$ and $\frac{|\nabla u|^2}{u}$  we see that
\begin{equation*}
\int_M|\Delta u|dg(t) \leq \int_M \Delta u dg(t) +C =C <\infty
\end{equation*}
for each $t \in [0,T)$, where the constant $C$  depends on $t$, but is uniform in each closed subinterval $[0,T']\subset [0,T)$; compare the proof of Lemma 7.2 in \cite{CTY}. It follows that
\begin{equation}
\int_M |P(u)| dg(t) \leq C_1 < \infty
\end{equation}
for each $t \in [0,T)$, where the constant $C_1$  depends on $t$, but is uniform in each closed subinterval $[0,T']\subset [0,T)$.

 On pp. 22-23 of \cite{Ku1}, by using  a family of cutoff functions $\phi_k$  constructed on pp. 17-18 of \cite{Ku1}  Kuang shows that
\begin{equation}
\frac{d}{dt}W(g(t),v(\cdot,t), \tau(t))=\frac{d}{dt}\int_M P(u)dg(t)
\end{equation}
when the right hand side makes sense. But note that the second equality in (3.25) on p. 23 of \cite{Ku1} needs a justification.  Here is a way to bypass it; compare \cite{Z12b}. From (4.2) we have
\begin{equation*}
\frac{d}{dt}\int_M P(u)\phi_kdg(t)=\int_M 2\tau  |Ric-Hess \ln u-\frac{g}{2\tau }|^2u\phi_kdg(t)-\int_MP(u)\Delta \phi_kdg(t),
 \end{equation*}
where $\phi_k$ is as  mentioned above, and
\begin{equation}
\begin{split}
& \int_{t_1}^{t_2}\int_M 2\tau  |Ric-Hess \ln u-\frac{g}{2\tau }|^2u\phi_kdg(t)dt\\
= & \int_M P(u)\phi_kdg(t_2)-\int_M P(u)\phi_kdg(t_1)+\int_{t_1}^{t_2}\int_M P(u)\Delta \phi_kdg(t)dt
\end{split}
\end{equation}
for $t_1,  t_2 \in [0,T)$.    As in Step 3 in the proof of Corollary 4.1 in \cite{Z12b},  with the aid of  (4.3) above and the property of $\phi_k$ on  p.18  in \cite{Ku1} and using Lebesgue's  dominated convergence theorem (for the RHS of (4.5)) and monotone convergence theorem (for the LHS of (4.5)), we can take limit in (4.5) as $k\rightarrow \infty$ and get
\begin{equation*}
 \int_{t_1}^{t_2}\int_M 2\tau  |Ric-Hess \ln u-\frac{g}{2\tau }|^2udg(t)dt\\
=  \int_M P(u)dg(t_2)-\int_M P(u)dg(t_1).
\end{equation*}
It follows that
\begin{equation*}
\frac{d}{dt}\int_M P(u)dg(t)=\int_M 2\tau  |Ric-Hess \ln u-\frac{g}{2\tau }|^2udg(t)
\end{equation*}
for any $t\in [0,T)$. Combining this with (4.4) we get (4.1).

The monotonicity of the $\mu$-functional is stated on p. 1847 in \cite{L}, but the condition under which it holds is not stated explicitly there. (Note that \cite{L}  cites Theorem 7.1 (ii) in \cite{CTY} which needs assumption $(a1)$ there.) With the additional condition on the injectivity radius we know that $\mu (g(t_2), \tau(t_2))$ is finite. By taking a minimizing sequence   of   the functional  $W(g(t_2), \cdot, T-t_2) $ on the set $\{  v\in C_0^\infty(M)   \hspace{2mm} | \hspace{2mm}  ||v||_{L^2(M,g(t_2))}=1\}$ and evolving  the square of each of its elements backward under the conjugate heat equation coupled to the Ricci flow, one can derive the monotonicity of the $\mu$-functional from that of the $W$-entropy (i.e. the formula (4.1))
as in the compact case.  \hfill{$\Box$}

\vspace*{0.4cm}

Using Proposition 4.1 one can extend  the uniform Sobolev inequality  along the Ricci flow on compact manifolds proved by Zhang \cite{Z07} (see also Ye \cite{Y}) to the noncompact case  as  in \cite{Ku1} and \cite{Ku2}. Note that on p. 36 of \cite{Ku1} Kuang used the method of differentiation under the integral sign; in the noncompact case this method needs a justification, but this can be done as in the proof of Lemma 2.2.2 in \cite{D}. By the way, \cite{Ku2} is not available to me; but see the reviews in MathSciNet and zbMATH.  Note that on p. 31 of \cite{Ku1} Kuang used the minimizer of the $W$-entropy to derive the monotonicity of the $\mu$-functional. As pointed out in Zhang \cite{Z12b} the minimizer of the $W$-entropy on a noncompact manifold does not always exist. However the monotonicity of the $\mu$-functional in the  situation in \cite{Ku1} holds true; see the proof of Proposition 4.1 above.

\vspace*{0.4cm}

Combining this uniform Sobolev inequality  along the Ricci flow on noncompact manifolds (and adapting Step 3 in the proof of Theorem 20 in \cite{Ku1}) and results in Section 3 one can extend some results in Zhang \cite{Z12a} and Bamler-Zhang \cite{BZ} including distance distortion estimates, construction of a cutoff function, heat kernel estimates,  a backward pseudolocality theorem  and a strong $\varepsilon$-regularity theorem for Ricci flow on compact manifolds to the following situation: $(M, (g(t))_{t \in [0,T)})$, $T< \infty$ being a complete solution to the Ricci flow   with  $\sup_{M\times [0,T']}|Rm|< \infty$  for any $0< T' < T$, $\sup_M|\nabla Rm_{g_0}|< \infty$ and $\sup_M|\nabla^2 Rm_{g_0}|< \infty$, where  the injectivity radius of  the initial metric $g_0$ is bounded away from
0. ( Note that on p. 411 of \cite{BZ}, the equality $\int \Delta K(\cdot, t)dg_t=0$ for the heat kernel is used. In the noncompact case, this needs justification; but this can be done by adapting the argument used by Kuang in deriving  (3.18) in \cite{Ku1}.)

\section{Appendix: Hamilton's comparison function and Shi's local derivative estimates}

In this appendix we recall Hamilton's construction of a comparison function and the application to Shi's local derivative estimates, see  \cite{H13}. We will clarify some points in \cite{H13} and add some details at certain places.

Fix $T >0$. Let $M$ be  a manifold  (without boundary) of dimension $n$.  Suppose $g(t)$ is a solution (not necessarily complete) to the Ricci flow  on $M \times [0, T]$.
 Fix $x_0 \in M$. Let $s=s(p,t):=d_t(x_0,p)$ be the  distance function between $x_0$ and $p\in M$ w.r.t. $g(t)$.  Choose a spacetime point $(p_0,t_0)\in M\times (0,T]$ with $p_0\neq x_0$. Assume that $\gamma_0$ is a minimal geodesic w.r.t. the metric $g(t_0)$ from  $x_0$ to $p_0$ parametrized by the arc length $\sigma$ (also w.r.t.  $g(t_0)$).

 First we want to construct a regularization of the function $s$ at $(p_0,t_0)$ if  $s$ is not smooth at $(p_0,t_0)$.  Let $u=u(\sigma)$, $\sigma \in [0,s(p_0,t_0)]$, be a smooth function with $|u|\leq 1$, $u(0)=0$, and $u(s(p_0,t_0))=1$. Given any point $p$ near $p_0$, there is a unique vector $\tilde{v}\in T_{p_0}M$ with $p=\exp_{p_0}^{g(t_0)}\tilde{v}$, where $\exp_{p_0}^{g(t_0)}$ is the exponential map at $p_0$ w.r.t. $g(t_0)$. Parallel translate (under the Levi-Civita connection of $g(t_0)$) $\tilde{v}$ along $\gamma_0$ back to $x_0$; the vector at $\gamma_0(\sigma)$ that we get in this process  will be denoted by $v(\sigma)$.  Let $\gamma_p^\dag$ be the curve from $x_0$ to $p$ defined by
\begin{equation*}
\gamma_p^\dag(\sigma):=\exp_{\gamma_0(\sigma)}^{g(t_0)}u(\sigma)v(\sigma),  \hspace{4mm} \sigma \in [0,s(p_0,t_0)].
\end{equation*}
For $t$ near $t_0$ let $\hat{s}(p,t)$ be the length of the curve $\gamma_p^\dag$ w.r.t. the metric $g(t)$.
Then by definition $\hat{s}$ is a smooth function defined for $(p,t)$ near $(p_0,t_0)$, and satisfies  $\hat{s}(p_0,t_0)=s(p_0,t_0)$ and  $\hat{s}(p,t)\geq s(p,t)$ for $(p,t)$ near $(p_0,t_0)$.

We claim that $|\nabla \hat{s}| =1$ at $(p_0,t_0)$.  To see this, we compute $\tilde{v}(\hat{s}(\cdot,t_0))$  for any vector $\tilde{v}\in T_{p_0}M$. As before, parallel translate (under the Levi-Civita connection of $g(t_0)$) $\tilde{v}$ along $\gamma_0$ back to $x_0$, and denote the vector at $\gamma_0(\sigma)$ that we get in this process  by $v(\sigma)$.  For any $\tau \in \mathbb{R}$ with $|\tau|$ sufficiently small, consider the curve $\gamma_\tau$ from  $x_0$ to $\exp_{p_0}^{g(t_0)}\tau \tilde{v}$ defined by
\begin{equation*}
\gamma_\tau(\sigma):=\exp_{\gamma_0(\sigma)}^{g(t_0)}\tau u(\sigma)v(\sigma),     \hspace{4mm} \sigma \in [0,s(p_0,t_0)].
\end{equation*}
 Denote the length of the curve $\gamma_\tau$ w.r.t. $g(t_0)$ by $L(\tau)$.  By definition
 \begin{equation*}
 L(\tau)=\hat{s}(\exp_{p_0}^{g(t_0)}\tau \tilde{v},t_0).
 \end{equation*}
 By the first variation formula for arc length and using the assumption that $\gamma_0$ is a geodesic w.r.t. $g(t_0)$ we have
  \begin{equation}
  \tilde{v}(\hat{s}(\cdot,t_0))=L'(0)=g_{t_0}(\gamma_0'(s(p_0,t_0)),\tilde{v}).
  \end{equation}
  So
  \begin{equation*}
  |\tilde{v}(\hat{s}(\cdot,t_0))| \leq |\tilde{v}|_{g(t_0)},
  \end{equation*}
 and
  \begin{equation}
 |\nabla \hat{s}| \leq 1
  \end{equation}
 at $(p_0,t_0)$, since $\tilde{v}\in T_{p_0}M$ is arbitrarily chosen.
 On the other hand, if we choose $\tilde{v}= \gamma_0'(s(p_0,t_0))$  in  (5.1), we get
 \begin{equation}
  \gamma_0'(s(p_0,t_0))(\hat{s}(\cdot,t_0))=1.
  \end{equation}
  It follows from (5.2) and (5.3) that
   \begin{equation}
  |\nabla \hat{s}|=1
   \end{equation}
   at $(p_0,t_0)$.

The following result of Hamilton \cite{H13} is complementary to Lemma 8.3 (a) in Perelman \cite{P}.

\begin{thm} \label{thm 5.1} (Hamilton \cite{H13})  Fix $T >0$. Let $M$ be  a manifold  (without boundary) of dimension $n$.  Suppose $g(t)$ is a solution (not necessarily complete) to the Ricci flow  on $M \times [0, T]$.
 Fix $(x_0,t_0)\in M\times (0,T]$ and $r>0$. Assume that  the closure of the open metric ball $B(x_0, t_0, \frac{\pi}{2}r)$ is compact, and $Ric_{g(t_0)}\leq \frac{n-1}{r^2}$ on $B(x_0, t_0, \frac{\pi}{2}r)$. Let  $s=s(p,t)$ be the distance function between $x_0$ and $p\in M$ w.r.t. $g(t)$.
 Then
  \begin{equation}
  \frac{\partial s}{\partial t} \geq \Delta s-\frac{n-1}{r}\varphi(\frac{s}{r})
 \end{equation}
at $t=t_0$ on $B(x_0, t_0, \frac{\pi}{2}r)\setminus \{x_0\}$  in the constructive comparison sense, where
\begin{equation*}
\varphi(y)=y+\cot y.
\end{equation*}
 That is, for any $p_0 \in B(x_0, t_0, \frac{\pi}{2}r)\setminus \{x_0\}$, one can construct a smooth function $\hat{s}$ in a spacetime neighborhood $U$ of $(p_0,t_0)$ with $\hat{s}\geq s$ in $U$,  $\hat{s}= s $ at $(p_0,t_0)$, $|\nabla \hat{s}|= 1$ at $(p_0,t_0)$, and
 \begin{equation}
  \frac{\partial \hat{s}}{\partial t} \geq \Delta \hat{s}-\frac{n-1}{r}\varphi(\frac{\hat{s}}{r})
 \end{equation}
 at $(p_0,t_0)$ in the classical sense.  Actually if the distance function $s$ is smooth at $(p_0,t_0)$,  the inequality (5.5) holds at $(p_0,t_0)$ in the classical sense.
   \end{thm}

  \noindent{\bf Proof} (Hamilton \cite{H13}). Fix $p_0 \in B(x_0, t_0, \frac{\pi}{2}r)\setminus \{x_0\}$.  We may assume that the function $s$ is not smooth at $(p_0,t_0)$; otherwise we just take $\hat{s}=s$, that is, the inequality (5.5) holds in the classical sense at a point $(p_0,t_0)$ (with $p_0 \in B(x_0, t_0, \frac{\pi}{2}r)\setminus \{x_0\}$) where $s$ is smooth, which can be proved by almost the same argument as below (cf. the proof of Lemma 8.3 (a) in Perelman \cite{P}).
   Let $\hat{s}$ be a function as constructed above (around $(p_0,t_0)$), which depends on another function $u$. Recall that  $\hat{s}$  is smooth in a spacetime neighborhood $U$ of $(p_0,t_0)$ with $\hat{s}\geq s$ in $U$,  $\hat{s}= s $ at $(p_0,t_0)$, and  $|\nabla \hat{s}|= 1$ at $(p_0,t_0)$  by (5.4).

    We have (see for example Section 12 of \cite{H95})
 \begin{equation}
 \frac{\partial \hat{s}}{\partial t}(p_0, t_0)=-\int_0^{s(p_0,t_0)}Ric_{g(t_0)}(\gamma_0'(\sigma),\gamma_0'(\sigma))d\sigma.
 \end{equation}
  Moreover by using the Fermi coordinates (w.r.t. $g(t_0)$) along the geodesic $\gamma_0$ and the second variation formula for arc length, we can compute
  \begin{equation}
 \Delta \hat{s}(p_0, t_0)=\int_0^{s(p_0,t_0)} ((n-1)u'(\sigma)^2-Ric_{g(t_0)}(\gamma_0'(\sigma),\gamma_0'(\sigma))u(\sigma)^2) d\sigma,
 \end{equation}
  cf. also for example Section 1 in Chapter I of \cite{SY}  and the proof of Lemma 8.3 (a) in Perelman \cite{P}.

   By (5.7) and (5.8) we have
   \begin{equation*}
\begin{split}
 & (\frac{\partial \hat{s} }{\partial t}-\Delta \hat{s}) (p_0,t_0) \\
= & -\int_0^{s(p_0,t_0)} ((n-1)u'(\sigma)^2+Ric_{g(t_0)}(\gamma_0'(\sigma),\gamma_0'(\sigma))(1-u(\sigma)^2)) d\sigma\\
\geq & -(n-1)\int_0^{s(p_0,t_0)} (u'(\sigma)^2+\frac{1}{r^2}(1-u(\sigma)^2)) d\sigma,
\end{split}
\end{equation*}
 where in the last inequality we use our  assumption on the Ricci curvature upper bound and the assumption $|u|\leq 1$.

The  functional
\begin{equation*}
 J=J(u):=\int_0^{s(p_0,t_0)} (u'(\sigma)^2+\frac{1}{r^2}(1-u(\sigma)^2)) d\sigma
 \end{equation*}
for smooth functions  $u$  on $[0, s(p_0,t_0)]$ with constraints $|u|\leq 1$,  $u(0)=0$, and $u(s(p_0,t_0))=1$, has a  minimizer,  denoted by $\hat{u}$, which solves   the Euler-Lagrange equation
\begin{equation*}
\frac{d^2u}{d\sigma^2}+\frac{1}{r^2}u=0
 \end{equation*}
with the same constraints.  So
 \begin{equation*}
\hat{u}= \hat{u}(\sigma)=\sin \frac{\sigma}{r}/ \sin \frac{s(p_0,t_0)}{r},    \hspace{4mm} \sigma \in [0,s(p_0,t_0)].
 \end{equation*}
Now we specify $\hat{s}$ by choosing $u$ to be $\hat{u}$.  Then the desired inequality (5.6) at $(p_0,t_0)$ follows  for this $\hat{s}$.
(We may  call this $\hat{s}$ the (Hamilton) regularization of $s$ at $(p_0,t_0)$.)

\hfill{$\Box$}

\vspace*{0.4cm}

\noindent {\bf Remark}.  As in  Lemma 8.3 (a) in Perelman \cite{P}, the heat inequality (5.5) for the distance function under the Ricci flow can also be understood in the barrier sense  (see for example \cite{C58}  and \cite{Do}), whose definition is more general.  Also note that the assumption that  the closure of the open metric ball $B(x_0, t_0, \frac{\pi}{2}r)$ is compact guarantees  the existence of a minimal geodesic w.r.t. the metric $g(t_0)$ connecting $x_0$ and $p_0$, so the $\hat{s}$ above can be constructed.

\vspace*{0.4cm}

\begin{cor} \label{cor 5.2} (Hamilton \cite{H13})  Fix $T >0$. Let $M$ be  a manifold  (without boundary) of dimension $n$.  Suppose $g(t)$ is a solution (not necessarily complete) to the Ricci flow  on $M \times [0, T]$.
 Fix $(x_0,t_0)\in M\times (0,T]$ and $r>0$. Assume that  the closure of the open metric ball $B(x_0, t_0, r)$ is compact, and $Ric_{g(t_0)}\leq \frac{n-1}{r^2}$ on $B(x_0, t_0, r)$. Then there exists a positive constant $A$ depending only on the dimension with the following property. Let  $s=s(p,t)$ be the distance function between $x_0$ and $p\in M$ w.r.t. $g(t)$, and
 \begin{equation*}
  \Phi=\frac{Ar^2}{(r^2-s^2)^2}+\frac{1}{t}
 \end{equation*}
 for $0 \leq s < r$ and $t$ near $t_0$.
 Then
  \begin{equation*}
  \frac{\partial \Phi}{\partial t} > \Delta \Phi-\Phi^2
 \end{equation*}
at $t=t_0$ on $B(x_0, t_0, r)$  in the constructive comparison sense.
 That is, for any $p_0 \in B(x_0, t_0, r)$, setting
 \begin{equation*}
  \hat{\Phi}=\frac{Ar^2}{(r^2-\hat{s}^2)^2}+\frac{1}{t}
 \end{equation*}
 for $(p,t)$ near $(p_0,t_0)$, where $\hat{s}=\hat{s}(p,t)$ is the regularization of $s$ at $(p_0,t_0)$  as  constructed in the proof of Theorem 5.1 in the case $p_0\neq x_0$  and $s$ is not smooth at $(p_0,t_0)$, and $\hat{s}=s$ otherwise, we have  $\hat{\Phi}\geq \Phi$ in a spacetime neighborhood of $(p_0,t_0)$, $\hat{\Phi}= \Phi$
 at $(p_0,t_0)$, and
  \begin{equation*}
  \frac{\partial \hat{\Phi}}{\partial t} > \Delta \hat{\Phi}-\hat{\Phi}^2
 \end{equation*}
 at $(p_0,t_0)$ in the classical sense.
   \end{cor}
\noindent{\bf Proof} (Hamilton \cite{H13}).  Set $\rho=\frac{2}{\pi}r$. Then
   \begin{equation*}
   Ric_{g(t_0)}\leq \frac{n-1}{\rho^2}  \hspace{4mm}  \text{on}  \hspace{4mm}  B(x_0,t_0,\frac{\pi}{2}\rho).
    \end{equation*}
    Given $p_0 \in B(x_0, t_0, r)= B(x_0,t_0,\frac{\pi}{2}\rho)$, if  $p_0\neq x_0$ and $s$ is not smooth at $(p_0,t_0)$ let $\hat{s}$ be the  regularization of $s$ at $(p_0,t_0)$  as constructed  in the proof of Theorem 5.1, otherwise let $\hat{s}=s$.
   By Theorem 5.1, if $p_0\neq x_0$ we have
\begin{equation*}
  \frac{\partial \hat{s}}{\partial t} \geq \Delta \hat{s}-\frac{n-1}{\rho}\varphi(\frac{\hat{s}}{\rho})
 \end{equation*}
 at $(p_0,t_0)$ in the classical sense.  Note that
  \begin{equation*}
\varphi(y) < y+\frac{1}{y}   \hspace{4mm}  \text{for}  \hspace{4mm} 0<y < \frac{\pi}{2}.
\end{equation*}
 It follows that if $p_0\neq x_0$ we have
\begin{equation}
  \frac{\partial \hat{s}}{\partial t} \geq \Delta \hat{s}-C(\frac{1}{\hat{s}}+\frac{1}{r})
 \end{equation}
 at $(p_0,t_0)$ in the classical sense, where $C$ is a positive constant depending only on the dimension.

  For a smooth even function $f=f(y)$ of one real variable $y$ defined in a neighborhood of $y=0$ which is  increasing for $y\geq 0$ we have
   \begin{equation*}
 \Delta f(\hat{s})=f'(\hat{s})\Delta \hat{s}+f''(\hat{s})|\nabla \hat{s}|^2 = f'(\hat{s})\Delta \hat{s}+f''(\hat{s})
 \end{equation*}
  at $(p_0,t_0)$  in the classical sense  when $p_0\neq x_0$, where we use the fact  $|\nabla \hat{s}|= 1$  at $(p_0,t_0)$ when $p_0\neq x_0$. Combining this with (5.9) and using the assumption $f'(y) \geq 0$ for $y\geq 0$ we see that
 \begin{equation}
 (\frac{\partial}{\partial t}-\Delta) f(\hat{s})\geq -f''(\hat{s})-C(\frac{1}{\hat{s}}+\frac{1}{r})f'(\hat{s})
 \end{equation}
 at $(p_0,t_0)$ in the classical sense  when $p_0\neq x_0$. Recall that by \cite{W} a smooth even function $f=f(y)$ of one real variable $y$ defined in a neighborhood of $y=0$ is in fact a smooth function of $y^2$, and the function $\frac{f'(y)}{y}$  extends smoothly over the point $y=0$.  Observe that the distance function $s$ is smooth at $(p_0,t_0)$ if $p_0\in B(x_0,t_0, \varepsilon)\setminus \{x_0\}$, where $\varepsilon $ is a sufficiently small positive number, so we can always take $\hat{s}=s$ for such $(p_0,t_0)$. Note also that the function $s^2$ is smooth at  $(x_0,t_0)$. So by continuity (5.10) also holds at  $(x_0, t_0)$  in the classical sense.

 Then by using (5.10) with
  \begin{equation*}
 f=f(y):=\frac{Ar^2}{(r^2-y^2)^2},   \hspace{4mm}  y\in (-r,r),
  \end{equation*}
 where $A$ is a positive constant,  a direct computation   shows that we can choose the  constant $A$ to depend only on the dimension, such that the function
 \begin{equation*}
 \hat{\Psi}:=\frac{Ar^2}{(r^2-\hat{s}^2)^2}
 \end{equation*}
 defined for $(p,t)$ near $(p_0,t_0)$   satisfies the strict inequality
\begin{equation*}
  \frac{\partial \hat{\Psi}}{\partial t} > \Delta \hat{\Psi}-\hat{\Psi}^2
 \end{equation*}
 at $(p_0,t_0)$ in the classical sense; here $p_0=x_0$ is allowed.  Hence the function
 \begin{equation*}
  \hat{\Phi}:=\frac{Ar^2}{(r^2-\hat{s}^2)^2}+\frac{1}{t}
 \end{equation*}
defined for $(p,t)$ near $(p_0,t_0)$ also satisfies the strict inequality
\begin{equation*}
  \frac{\partial \hat{\Phi}}{\partial t} > \Delta \hat{\Phi}-\hat{\Phi}^2
 \end{equation*}
 at $(p_0,t_0)$ in the classical sense.  Clearly we also have  $\hat{\Phi}\geq \Phi$ in a spacetime neighborhood of $(p_0,t_0)$, and $\hat{\Phi}= \Phi$
 at $(p_0,t_0)$.
   \hfill{$\Box$}

\vspace*{0.4cm}

\noindent {\bf Remark}. In the static metric case a  function  (of the same form as $\Phi$ above)
\begin{equation*}
   h= \frac{Ar^2}{(r^2-s^2)^2} + \frac{1}{t}
\end{equation*}
(using our notation) appears in the proof of Theorem F1.4 in Hamilton \cite{H97} which satisfies
\begin{equation*}
  \frac{\partial h}{\partial t} \geq \Delta h-h^2
\end{equation*}
in a small ball of  radius $r$ with some suitably chosen positive constant $A$.

\vspace*{0.4cm}

\begin{thm} \label{thm 5.3} (Shi's local derivative estimates, see \cite{H13})  Fix $T >0$. Let $M$ be  a manifold  (without boundary) of dimension $n$. Suppose $g(t)$ is a solution (not necessarily complete) to the Ricci flow  on $M \times [0, T]$.  Fix $x_0 \in M$ and $r>0$.  Assume that  the parabolic cylinder $PB_r(x_0,T)$ is compact, and $|Rm|\leq K$ on  $PB_r(x_0,  T)$ with (the  constant) $K\leq \frac{1}{r^2}$. Then for any $i\geq 1$,
\begin{equation*}
|\nabla^i Rm| \leq C_i K(\frac{1}{r^{i}}+\frac{1}{t^{i/2}})    \hspace{2mm}  \text{on}   \hspace{2mm}    PB_\frac{r}{2^i}(x_0,  T)\setminus \{(x, 0) | x \in M\},
\end{equation*}
where the constant $C_i$ depends only on $i$ and the dimension.
\end{thm}

\noindent{\bf Proof}.  We do induction on $i$. For the case $i=1$ we reproduce  Hamilton's argument in  \cite{H13}. Let $\alpha$ be a large positive constant, which will be chosen later, and let
\begin{equation*}
b=\frac{1}{(\alpha+1)^2K^4}  \hspace{4mm}  \text{and}   \hspace{4mm}  B=\frac{K}{\sqrt{\alpha}}.
\end{equation*}
\vspace*{0.4cm}
Set
\begin{equation*}
F=b(\alpha K^2+|Rm|^2)|\nabla Rm|^2-B.
\end{equation*}
Then by  choosing the constant $\alpha$ large enough (depending only on the dimension) we get
\begin{equation*}
\frac{\partial}{\partial t}F \leq \Delta F -F^2.
\end{equation*}
We claim that
\begin{equation*}
F < \Phi  \hspace{2mm}  \text{for}   \hspace{2mm}    x\in B(x_0,t,r)  \hspace{2mm}  \text{and}   \hspace{2mm}  t\in (0,T],
\end{equation*}
where $\Phi$ is the comparison function in Corollary 5.2.  The proof uses a maximum principle argument, and goes as follows. Note that $\Phi \rightarrow \infty$ as $(p,t)$ tends to the parabolic boundary of $PB_r(x_0,T)$. But $F$ is bounded on $PB_r(x_0,T)$. So the set
\begin{equation*}
D:=\{(x,t)| x\in B(x_0,t,r), t\in (0,T], F(x,t)\geq \Phi(x,t)\}
 \end{equation*}
 is a compact subset of the set $\{(x,t)| x\in B(x_0,t,r), t\in (0,T]\}$. (Here we also use the assumption  that the parabolic cylinder $PB_r(x_0,T)$ is compact.) It turns out  that the set $D$ is empty. Argue by contradiction.  Otherwise we can choose a point $(p_0,t_0) \in D$ with $t_0$ minimal.  Then $F\leq \Phi$ when $t\leq t_0$, and $F=\Phi$ at $(p_0,t_0)$. We can construct a function $\hat{\Phi}$ as in Corollary 5.2 with $\hat{\Phi}\geq \Phi$ in a spacetime neighborhood of $(p_0,t_0)$, $\hat{\Phi}= \Phi$ at $(p_0,t_0)$, and
 \begin{equation*}
  \frac{\partial \hat{\Phi}}{\partial t} > \Delta \hat{\Phi}-\hat{\Phi}^2
 \end{equation*}
 at $(p_0,t_0)$ in the classical sense.  On the other hand, we still have $F\leq \hat{\Phi}$ when $t\leq t_0$, and $F=\hat{\Phi}$ at $(p_0,t_0)$ (in particular, $(p_0,t_0)$ is a maximum point of the function $F-\hat{\Phi}$ restricted to the time slice $t=t_0$). So we have
 \begin{equation*}
  \frac{\partial F}{\partial t} \geq \frac{\partial \hat{\Phi}}{\partial t}  \hspace{4mm}   \text{and}  \hspace{4mm} \Delta F \leq \Delta \hat{\Phi}
 \end{equation*}
 at $(p_0,t_0)$. It follows that
  \begin{equation*}
  \frac{\partial F}{\partial t} \leq \Delta F-F^2 \leq \Delta \hat{\Phi}-\hat{\Phi}^2 < \frac{\partial \hat{\Phi}}{\partial t} \leq \frac{\partial F}{\partial t}
 \end{equation*}
  at $(p_0,t_0)$.  A contradiction.  So the set $D$ is empty, and the claim holds.  The desired conclusion in the case $i=1$ follows easily.

 Now suppose we have
 \begin{equation}
|\nabla^i Rm| \leq C_i K(\frac{1}{r^{i}}+\frac{1}{t^{i/2}})    \hspace{4mm}  \text{on}   \hspace{2mm}    PB_\frac{r}{2^i}(x_0,  T)\setminus \{(x, 0) | x \in M\},
\end{equation}
for $1\leq i \leq k$, where the constant $C_i$ depends only on $i$ and the dimension. For the case $i=k+1$,
 let
 \begin{equation*}
S_{k+1}=(B_{k+1}K^2(\frac{1}{r^{2k}}+\frac{1}{t^k})+|\nabla^k Rm|^2)|\nabla^{k+1} Rm|^2,
\end{equation*}
 where $B_{k+1}$ is a  constant to be chosen. By choosing $B_{k+1}$ large enough (depending only on $k$ and the dimension) and using (5.11) we have
 \begin{equation*}
\begin{split}
&(\frac{\partial}{\partial t}-\Delta)S_{k+1} \\
\leq &- \frac{S_{k+1}^2}{(B_{k+1}+1)^2K^4(\frac{1}{r^{2k}}+\frac{1}{t^k})^2}+ CB_{k+1}^2K^5(\frac{1}{r^{2(2k+1)}}+\frac{1}{t^{2k+1}})
\end{split}
\end{equation*}
 on  $\{(x,t)| x\in B(x_0, t, \frac{r}{2^k}), t\in (0, T]\}$,    where the constant $C$ depends only on $k$ and the dimension;  see p. 197 in \cite{CZ}.  (By the way, it seems that  there are two  typos in line 2 from the bottom on p. 197 in \cite{CZ}.)

 Now  setting $F_{k+1}=\beta S_{k+1}/v^k$, where  $\beta=\frac{1}{2(C+2k^2)(B_{k+1}+1)^2K^4}$, and $v=\frac{1}{r^2}+\frac{1}{t}$, we get
\begin{equation*}
(\frac{\partial}{\partial t}-\Delta)F_{k+1}\leq -\frac{F_{k+1}^2}{v^k}+v^{k+2}
\end{equation*}
on  $\{(x,t)| x\in B(x_0, t, \frac{r}{2^k}), t\in (0, T]\}$; see p. 198 in \cite{CZ}.

On the other hand, exactly as in the proof of Theorem 2.1,  using the comparison function $\Phi$ in Corollary 5.2 we can construct a function $\Phi_{k+1}$ satisfying

 \begin{equation*}
(\frac{\partial}{\partial t}-\Delta)\Phi_{k+1} > -\frac{\Phi_{k+1}^2}{v^k}+v^{k+2}
\end{equation*}
  on  $\{(x,t)| x\in B(x_0, t, \frac{r}{2^k}), t\in (0, T]\}$ in the constructive comparison sense.

Observe  that  $\Phi_{k+1} \rightarrow \infty$ as $(x,t)$ tends to the parabolic boundary of $PB_{\frac{r}{2^k}}(x_0,T)$, but $F_{k+1}$ is bounded  on   $\{(x,t)| x\in B(x_0, t, \frac{r}{2^k}), t\in (0, T]\}$  (note that $\frac{1}{r^{2k}}+\frac{1}{t^k} < v^k$). So near the parabolic boundary of $PB_{\frac{r}{2^k}}(x_0, T)$ we have
 $F_{k+1} < \Phi_{k+1}$.
 Using a maximum principle argument as in the case $i=1$ above we get that $F_{k+1} < \Phi_{k+1}$ everywhere on
$\{(x,t)| x\in B(x_0, t, \frac{r}{2^k}), t\in (0, T]\}$. Then the desired conclusion in the case $i=k+1$ follows easily.
 \hfill{$\Box$}

\vspace*{0.4cm}

\begin{thm} \label{thm 5.4} (Shi's local derivative estimates with initial derivative bounds, see \cite{H13}) If in addition to the assumptions in Theorem 5.3  we also assume
\begin{equation*}
|\nabla^i Rm| \leq  \frac{K}{r^{i}}    \hspace{4mm}  \text{on}   \hspace{2mm}    B(x_0,0,r),  \hspace{2mm}  1\leq i \leq p,
\end{equation*}
for some $p\geq 1$, then we have
\begin{equation*}
|\nabla^i Rm| \leq C_i \frac{K}{r^{i}}    \hspace{4mm}  \text{on}   \hspace{2mm}    PB_\frac{r}{2^i}(x_0,  T),  \hspace{2mm}  1\leq i \leq p,
\end{equation*}
where the constant $C_i$ depends only on $i$ and the dimension.
\end{thm}

\noindent{\bf Proof}. We only need to slightly modify  the proof of Theorem 5.3.  For example in the case $p=1$,  we use  the comparison function
$\Psi:=\frac{Ar^2}{(r^2-s^2)^2}$  instead of  $\Phi$.  The higher order case is similar: We use $(B_iK^2\frac{1}{r^{2(i-1)}}+|\nabla^{i-1} Rm|^2)|\nabla^i Rm|^2$ as $S_i$ above,  and  $\beta_{i} \Psi_{i}^{i}$ instead of $\Phi_{i}$ ($1\leq i \leq p$) as  the (space-time) comparison function; cf. the proof of Theorem 2.1 and the remark after it.   \hfill{$\Box$}

\vspace*{0.4cm}

\noindent {\bf Acknowledgements}.   I would like to thank Professor Qi S. Zhang for sending me a copy of \cite{Ku1} and bringing the paper \cite{RV} to my attention. I would also like to thank Dr. Yu Li for answering my question on his paper \cite{L} related to the monotonicity of the $\mu$-functional on noncompact manifolds. I'm also grateful to the referee for the comments and suggestions which help improve the presentation of the paper. I was partially supported by  Laboratory of Mathematics and Complex Systems, Ministry of Education,  and by  Beijing Natural Science Foundation (Z190003).

% ----------------------------------------------------------------

\hspace *{0.4cm}

\vspace *{0.4cm}

Laboratory of Mathematics and Complex Systems (Ministry of Education),

School of Mathematical Sciences, Beijing Normal University,

Beijing 100875, P.R. China

 E-mail address: hhuang@bnu.edu.cn

\end{document}